\documentclass{amsart}

\newcommand\cF{{\mathcal{F}}}
\newcommand\cG{{\mathcal{G}}}

\newcommand\cL{{\mathcal{L}}}

\newcommand\cS{{\mathcal{S}}}

\newcommand\ZZ{\mathbb{Z}}

\newcommand\CC{\mathbb{C}}

\newcommand\NN{\mathbb{N}}

\theoremstyle{plain}
\newtheorem{thm}{Theorem}[section]

\newtheorem{lem}[thm]{Lemma}

\newtheorem{prop}[thm]{Proposition}

\theoremstyle{definition}
\newtheorem{defn}[thm]{Definition}
\newtheorem{ex}[thm]{Example}
\newtheorem{hyp}[thm]{Hypotheses}

\theoremstyle{remark}
\newtheorem{rem}[thm]{Remark}

\begin{document}

\title[Noncommutative Hamiltonian dynamics on foliated manifolds]{Noncommutative Hamiltonian dynamics \\ on foliated manifolds}
\author{Yuri A. Kordyukov}
\address{Institute of Mathematics,\\ Russian Academy of Sciences,\\ 112 Chernyshevsky
str.\\ 450008 Ufa\\ Russia} \email{yurikor@matem.anrb.ru}

\thanks{Partially supported by the Russian Foundation of Basic
Research (grant 07-01-00081-a).}

\hyphenation{geo-de-sic}

\bibliographystyle{plain}

\subjclass{Primary 58B34; Secondary 37J05, 53C12}

\keywords{Noncommutative geometry, Poisson structures, Hamiltonian
flows, foliations}

\date{}

\dedicatory{Dedicated to Stephen Smale on his 80th birthday}

\begin{abstract}
First, we review the notion of a Poisson structure on a
noncommutative algebra due to Block-Getzler and Xu and introduce a
notion of a Hamiltonian vector field on a noncommutative Poisson
algebra. Then we describe a Poisson structure on a noncommutative
algebra associated with a transversely symplectic foliation and
construct a class of Hamiltonian vector fields associated with this
Poisson structure.
\end{abstract}
\maketitle


\section{Introduction}
The class of Hamiltonian flows on a symplectic manifold is an
important class of dynamical systems. It naturally arises as a
geometric manifestation of Hamilton's equations in classical
mechanics. Hamiltonian flows can be defined more generally on an
arbitrary Poisson manifold.

The purpose of this paper is to discuss the notion of a Hamiltonian
flow in the framework of Poisson geometry on a particular class of
singular symplectic manifolds, namely, on the leaf space of a
transversely symplectic compact foliated manifold. In this case, it
is natural to use the notions and methods of noncommutative
differential geometry initiated by Alain Connes \cite{Co}.

There are several fundamental ideas, which lie in the base of
noncommutative geometry. The first of them is to pass from geometric
spaces to algebras of functions on these spaces and translate basic
geometric and analytic notions and constructions into the algebraic
language. An application of this idea in Poisson geometry leads to a
notion of Poisson algebra.

The next idea is that, in many cases (especially, in those cases
when the classical commutative algebra of functions is small or has
a bad structure), it is useful to consider as its analogue some
noncommutative algebra. This gives rise to the need to extend the
basic geometric and analytic definitions to general noncommutative
algebras. A noncommutative analogue of a Poisson structure was
introduced independently by Block and Getzler \cite{Block-Ge} and Xu
\cite{Xu}. They used some ideas from deformation theory of
associative algebras.

In \cite{Co79}, Connes associated with an arbitrary foliated
manifold $(M,\cF)$ the $C^*$-algebra $C^*(M,\cF)$, which can be
naturally considered as a noncommutative analogue of the algebra of
continuous functions on the leaf space $M/\cF$ of the foliation.
When the foliation has a transverse symplectic structure, there is a
natural noncommutative Poisson structure defined on a dense
subalgebra of the $C^*$-algebra $C^*(M,\cF)$, which was constructed
in \cite{Block-Ge}.

In this paper, we introduce the notion of a Hamiltonian vector field
associated with a noncommutative Poisson structure and construct a
class of Hamiltonian vector fields on the $C^*$-algebra $C^*(M,\cF)$
associated with a transversally symplectic foliation $(M,\cF)$. This
study was partially motivated by our investigations of transversally
elliptic operators on foliated manifolds, related trace formulae and
the corresponding classical dynamics
\cite{dgtrace-eng,egorgeo,matrix-egorov}. In particular, it follows
from the results of the paper that the dynamical systems on
foliation algebras, which appear in the Egorov type theorems for
transversally elliptic operators proved in
\cite{egorgeo,matrix-egorov}, are Hamiltonian flows.

We remark that there is a notion of a noncommutative symplectic
manifold introduced by Kontsevich \cite{kontsevich93} and further
developed by Ginzburg in \cite{Ginzburg,Ginzburg:lectures}. One can
define a notion of a Hamiltonian vector field on a noncommutative
symplectic manifold. Moreover, one can show that a transversally
symplectic foliation gives rise to a noncommutative symplectic
manifold, and the noncommutative vector fields constructed in our
paper are Hamiltonian vector fields on this noncommutative
symplectic manifold. These results will be discussed elsewhere.

The paper is organized as follows. First, we review the notions of a
Poisson structure on an associative algebra and a Hamiltonian vector
field associated with a noncommutative Poisson structure. Next, we
describe the noncommutative geometry of the leaf space of a foliated
manifold and the noncommutative Poisson structure of a transversely
symplectic foliated manifolds and construct a class of
noncommutative Hamiltonian flows on these manifolds.

We refer the reader to \cite{Vaisman} for basic notions of Poisson
geometry and to the survey paper \cite{survey} for information and
references on noncommutative geometry of foliations.

\section{Preliminaries on noncommutative Poisson geometry}
In this Section, we review some basic notions related with
noncommutative Poisson structures on associative algebras, following
\cite{Block-Ge,Xu}.

\subsection{Noncommutative Poisson structures}
Let $A$ be an associative algebra over $\CC$. The space of
Hochschild $k$-cochains on $A$ is
$C^k(A,A)=\operatorname{Hom}(A^{\otimes k},A)$. The differential $b
: C^k(A,A)\to C^{k+1}(A,A)$ is given by
\begin{multline*}
(b c)(a_1,\cdots,a_{k+1})=a_1c(a_2,\cdots,a_{k+1})\\
+\sum_{i=1}^{k}(-1)^ic(a_1,\cdots,
a_ia_{i+1},\cdots,a_{k+1})+(-1)^kc(a_1,\cdots,a_k)a_{k+1}.
\end{multline*}
The cohomology of the complex $(C^*(A,A),b)$ is called the
Hochschild cohomology $H^*(A,A)$ of $A$. For example, $H^0(A,A)$ is
just the center of $A$, and $H^1(A,A)$ is the space
$\operatorname{Out}(A)=\operatorname{Der}(A)/\operatorname{Inn}(A)$
of outer derivations of $A$.

We`define a pre-Lie product on $C^\ast(A,A)$. For any $U\in
C^u(A,A)$ and $V\in C^v(A,A)$, $U\ast V\in C^{u+v-1}(A,A)$ is given
by
\begin{multline*}
(U\ast V)(a_1, \cdots , a_{u+v-1})\\ =\sum_{i=1}^u (-1)^{(i-1)(v-1)
}U(a_1, \cdots, a_{i-1}, V(a_{i},\cdots, a_{i+v-1}), a_{i+v},
\cdots, a_{u+v-1}).
\end{multline*}

The Gerstenhaber bracket \cite{Gersten63} is defined to be the
commutator of the pre-Lie bracket: for any $U\in C^u(A,A)$ and $V\in
C^v(A,A)$, $[U, V]\in C^{u+v-1}(A,A)$ is given by
\[
[U,V]=U\ast V-(-1)^{(u-1)(v-1)}V\ast U.
\]
The Gerstenhaber bracket is a generalization of the usual
Schouten-Nijenhuis brackets of multivector fields.

\begin{defn}\label{d:Poisson}
A Poisson structure on $A$ is a Hochschild two-cocycle $\Pi\in
Z^2(A,A)$ such that $[\Pi.\Pi]$ is a three-boundary, that is, $\Pi$
is a homomorphism $\Pi : A\otimes A\rightarrow A$ such that
\begin{enumerate}
\item for any $a_1,a_2,a_3\in A$,
\begin{equation}\label{P1}
a_1\Pi(a_2,a_3)-\Pi(a_1a_2,a_3)+\Pi(a_1,a_2a_3)-\Pi(a_1,a_2)a_3=0;
\end{equation}
($\Pi\in Z^2(A,A)$; an analogue of the Leibniz's rule).
\item there is a homomorphism
$\Pi_1:A\otimes A\rightarrow A$ such that, for any $a_1,a_2,a_3\in
A$,
\begin{multline}\label{P2}
\Pi(a_1,\Pi(a_2,a_3))-\Pi(\Pi(a_1,a_2),a_3) \\=
a_1\Pi_1(a_2,a_3)-\Pi_1(a_1a_2,a_3)+
\Pi_1(a_1,a_2a_3)-\Pi_1(a_1,a_2)a_3.
\end{multline}
($[\Pi.\Pi]$ is a three-boundary; an analogue of the Jacobi rule).
\end{enumerate}
\end{defn}

\begin{ex}
Let $M$ be a compact smooth manifold. Recall that any Poisson
bracket $\{\cdot,\cdot\}$ on $M$ is determined by a Poisson bivector
$\Lambda \in C^\infty(M,\Lambda^2TM)$:
\[
\{f,g\}=\langle \Lambda, df\wedge
dg\rangle=\sum_{ij}\Lambda^{ij}\frac{\partial f}{\partial
x^i}\frac{\partial g}{\partial x^j}, \quad f,g, \in C^\infty(M),
\]
where $\{x^i\}$ are local coordinates on $M$. In local coordinates,
$\Lambda$ has to satisfy the following condition
\[
\sum_\alpha\left(\Lambda^{\alpha i}\frac{\partial
\Lambda^{jk}}{\partial x^\alpha}+\Lambda^{\alpha j}\frac{\partial
\Lambda^{ki}}{\partial x^\alpha}+\Lambda^{\alpha k}\frac{\partial
\Lambda^{ij}}{\partial x^\alpha}\right)=0
\]
for any $i$, $j$ and $k$.

An invariant meaning of this identity is provided by the
Schouten-Nijenhuis bracket, which is a bilinear local type extension
of the Lie derivative $L_X$ to an operation
\[
[\cdot,\cdot] : C^\infty(M,\Lambda^pTM)\times
C^\infty(M,\Lambda^qTM) \to C^\infty(M,\Lambda^{p+q-1}TM).
\]
A bivector $\Lambda \in C^\infty(M,\Lambda^2TM)$ is a Poisson
bivector if and only if $[\Lambda,\Lambda]=0$.

Consider the commutative algebra $A=C^\infty(M)$ of smooth functions
on $M$. Then there is an isomorphism \cite{Vey75}
\[
H^*(C^\infty(M),C^\infty(M))=C^\infty(M,\Lambda^*TM).
\]
For any $p$-vector field $X_1\wedge\ldots\wedge X_p\in
C^\infty(M,\Lambda^pTM)$, the corresponding cocycle $c\in
H^p(C^\infty(M),C^\infty(M))$ is given by
\[
c(f_1,\ldots,f_p)=\det\|X_if_j\|_{i,j=1}^p, \quad f_1,\ldots,f_p\in
C^\infty_c(M).
\]

The bilinear map
\[
\Pi(f,g)=\{f,g\},\quad f,g\in C^\infty(M),
\]
is a Poisson structure on the algebra $C^\infty(M)$. The
corresponding homomorphism $\Pi_1$ is described as follows. Let
$\nabla$ be a torsion-free connection on $M$. It induces the
covariant derivative $\nabla : C^\infty(M,T^*M) \rightarrow
C^\infty_c(M,T^*M\otimes T^*M) $ on $T^*M$. We can consider the
composition of this operator with the de Rham differential $d :
C^\infty(M)\to C^\infty(M,T^*M)$:
\[
\nabla^2=\nabla\circ d : C^\infty(M)\to C^\infty(M,S^2T^*M).
\]
The operator $\nabla^2$ takes values in $C^\infty(M,S^2T^*M)$ since
$\nabla$ is torsion-free.

The corresponding two-cochain $\Pi_1$ is defined by
\[
\Pi_1(f,g)=\langle \Lambda\otimes \Lambda, \nabla^2f \otimes
\nabla^2g\rangle, \quad f,g, \in C^\infty(M),
\]
where the pairing of the tensor $\Lambda\otimes \Lambda$ with
$\alpha\otimes \beta \in S^2T^*M\otimes S^2T^*M$ is given by the
formula
\[
\langle \Lambda\otimes \Lambda, \alpha \otimes \beta\rangle
=\Lambda^{ij}\Lambda^{kl}\alpha_{ik}\beta_{jl}.
\]

For any $f\in C^\infty(M)$, the map $g\mapsto \{f,g\}$ is a
derivation of $C^\infty(M)$, Therefore, there exists a well defined
vector field $X_f$ on $M$ such that
\[
\{f,g\}=X_fg=-X_gf=dg(X_f)=-df(X_g).
\]
$X_f$ is called the Hamiltonian vector field of $f$.
\end{ex}

\begin{ex}
One of the basic examples in noncommutative differential geometry is
the noncommutative two-torus $A_\theta$. The algebra $A_\theta$ is
generated by two elements $U$ and $V$, satisfying the relation
\[
VU=e^{2\pi i \theta} U V.
\]
A generic element of $A_\theta$ can be represented as a formal power
series
\[
a=\sum_{(n,m)\in\ZZ^2}a_{nm}U^nV^m,
\]
where $a_{nm}\in \cS(\ZZ^2)$ is a rapidly decreasing sequence (that
is, for any natural $k$ we have
$\sup_{(n,m)\in\ZZ^2}(|n|+|m|)^k|a_{nm}|<\infty$).

$A_\theta$ is a locally convex topological algebra under the
topology generated by the seminorms
\[
p_k(a)=\sup_{(n,m)\in\ZZ^2}(|n|+|m|)^k|a_{nm}|, \quad k\in \NN.
\]

There are two canonical derivations $\delta_1$ and $\delta_2$ on
$A_\theta$ given, respectively, by
\[
\delta_1(U^nV^m)=2\pi i n U^nV^m, \quad \delta_2(U^nV^m)=2\pi i m
U^nV^m.
\]
It is easy to see that $[\delta_1,\delta_2]=0$.

By \cite[Theorem 4.1]{Xu}, there is a canonical Poisson structure
$\Pi\in Z^2(A_\theta,A_\theta)$ on $A_\theta$ defined by
\[
\Pi(a_1,a_2)=\delta_1(a_1)\delta_2(a_2), \quad a_1, a_2\in A_\theta.
\]
In particular, the Jacobi rule \eqref{P2} holds with
\[
\Pi_1(a_1,a_2)=-\frac12\delta^2_1(a_1)\delta^2_2(a_2), \quad a_1,
a_2\in A_\theta.
\]
\end{ex}

More examples of noncommutative Poisson structures can be found in
\cite{Block-Ge,Xu,tang:gafa,Halbout-Tang}.

\subsection{Hamiltonian dynamics}
For a given noncommutative Poisson structure $\Pi$ on an associative
algebra $A$ over $\CC$, we denote the center of $A$ by $C$, i.e.
$C=H^0(A,A)$.

\begin{defn}\label{d:ham}
For any element $c$ of $C$, the Hamiltonian derivation of $A$
associated to $c$ is defined as
\[
X_c=\frac12[\Pi,c]\in H^1(A,A),
\]
or equivalently
\[
X_c(a)=\frac12\left(\Pi(c,a)-\Pi(a,c)\right),\quad a\in A.
\]
\end{defn}

\begin{rem}
It is impossible, in general, to associate a Hamiltonian derivation
to an arbitrary element of $A$ due to the lack of outer derivations
in $A$.
\end{rem}

One can introduce a bracket on the center $C$ as follows: For any
$c$ and $e$ in $C$
\begin{equation}\label{e:bracket}
\{c,e\}=[X_c,e]\in H^0(A,A)=C.
\end{equation}

We have the following properties (see \cite[Proposition 2.1]{Xu}).

\begin{prop}
For any $c$ and $e$ in $C$
\begin{enumerate}
  \item $L_{X_c}\Pi=0$;
  \item $[X_c,X_e]=-X_{\{c,e\}}$;
  \item $C$ together with the bracket $\{\cdot,\cdot\}$ introduced
  above becomes a Poisson algebra in the usual sense.
\end{enumerate}
\end{prop}

\begin{rem}
It is easy to see that, for a compact Poisson manifold $M$, the
bracket on the commutative algebra $C^\infty(M)$ defined by
\eqref{e:bracket} coincides with the bracket on $C^\infty(M)$ given
by the Poisson structure, and the Hamiltonian derivation $X_f$ of
$C^\infty(M)$ associated to $f\in C^\infty(M)$ by
Definition~\ref{d:ham} is determined by the classical Hamiltonian
vector field with Hamiltonian $f$.
\end{rem}

\section{Transverse geometry of foliations}
Throughout in this Section, $(M,{\cF})$ is a compact foliated
manifold, $\operatorname{dim} M=n, \operatorname{dim} \cF=p, p+q=n$.
We will consider foliated charts $\phi: U\subset M\rightarrow
I^p\times I^q$ on $M$ with coordinates $(x,y)\in I^p\times I^q$ ($I$
is the open interval $(0,1)$) such that the restriction of $\cF$ to
$U$ is given by the level sets $y={\rm const}$. We will use the
following notation: $T\cF$ is the tangent bundle of $\cF$;
$\tau=TM/T\cF$ is the normal bundle of $\cF$; $N^*{\mathcal
F}=\{\xi\in T^*M : \langle \xi, X\rangle =0 \ \forall X\in T\cF\}$
is the conormal bundle of $\cF$.

\subsection{Transverse symplectic structures}
A transverse symplectic structure on a foliated manifold
$(M,{\mathcal F})$ is given by a covering $\{U_i,\phi_i\}$ by
foliated charts, $\phi_i:U_i\to I^p\times I^q$, and by a family of
symplectic forms $\omega_i$ on local bases $I^q$ such that for any
coordinate transformation
$$
\phi_{ij}(x,y)=(\alpha_{ij}(x,y), \gamma_{ij}(y)), \quad x\in I^p,
y\in I^q,
$$
the map $\gamma_{ij}$ preserves the symplectic structure,
$\omega_j=\gamma^*_{ij}\omega_i$.

A manifold $M$ is called presymplectic, if it is endowed with a
closed two-form $\omega$ of constant rank.

One can show \cite{Block-Ge} that presymplectic structures are
essentially the same as transverse symplectic structures. More
precisely, if $M$ is a presymplectic manifold and $F\subset TM$ is
the subbundle on which $\omega$ vanishes, then $F$ is integrable and
thus defines a foliation $\cF$ on $M$. The restrictions of $\omega$
to the local bases of foliated charts on $M$ define a transverse
symplectic structure on the foliated manifold $(M,{\mathcal F})$. On
the other hand, a transverse symplectic structure on a foliated
manifold $(M,{\mathcal F})$ determines in a unique manner a
presymplectic structure on $M$ such that $F=T\cF$ is the kernel of
$\omega$.

\subsection{Foliation algebras}
Here we will describe a noncommutative algebra, which can be
considered, according to noncommutative geometry, as an algebra of
functions on the leaf space $M/\cF$ of a foliation $(M,\cF)$. As a
vector space, this algebra is the space $C^\infty_c(G)$ of smooth
compactly supported functions on the holonomy groupoid $G$ of the
foliation. Therefore, we recall first the notion of holonomy
groupoid.

Consider the equivalence relation $\sim_h$ on the set of continuous
leafwise paths $\gamma:[0,1]\rightarrow M$, setting $\gamma_1\sim_h
\gamma_2$, if $\gamma_1$ and $\gamma_2$ have the same initial and
final points and the same holonomy maps: $h_{\gamma_1} =
h_{\gamma_2}$. The holonomy groupoid $G$ is the set of
$\sim_h$-equivalence classes of leafwise paths. The set of units
$G^{(0)}$ is  $M$. The multiplication in $G$ is given by the product
of paths. The corresponding range and source maps $s,r:G\rightarrow
M$ are given by $s(\gamma)=\gamma(0)$ and $r(\gamma)=\gamma(1)$.
Finally, the diagonal map $\Delta:M\rightarrow G$ takes any $x\in M$
to the element in $G$ given by the constant path $\gamma(t)=x, t\in
[0,1]$. To simplify the notation, we will identify $x\in M$ with
$\Delta(x)\in G$.

For any $x\in M$ the map $s$ takes the set $G^x=r^{-1}(x)$ onto the
leaf $L_x$ through $x$. The group $G^x_x=s^{-1}(x)\cap r^{-1}(x)$
coincides with the holonomy group of $L_x$. The map
$s:G^x\rightarrow L_x$ is the covering map associated with the group
$G^x_x$, called the holonomy covering.

The holonomy groupoid $G$ has a structure of a smooth (in general,
non-Hausdorff and non-paracompact) manifold of dimension $2p+q$
\cite{Co79}. A local coordinate system on $G$, denoted by
$W(\phi,\phi')$, is determined by a pair of compatible foliated
charts $\phi$ and $\phi'$ on $M$. The coordinates in $W(\phi,\phi')$
will denote by $(x,x',y)\in I^p\times I^p\times I^q$.

Let us fix a positive smooth leafwise density $\alpha \in
C^\infty(M,|T\cF|)$. For any $x\in M$, we define a positive Radon
measure $\nu ^{x}$ on $G^x$ to be the lift of the restriction of
$\alpha $ to $L_x$ by the holonomy cover $s:G^x\to L_x$. The
structure of an involutive algebra on $C^{\infty}_c(G)$ is defined
by
\begin{align*}
k_1\ast k_2(\gamma)&=\int_{G^x} k_1(\gamma_1)
k_2(\gamma^{-1}_1\gamma)\,d\nu^x(\gamma_1),\quad \gamma\in G^x,\\
k^*(\gamma)&=\overline{k(\gamma^{-1})}, \quad \gamma\in G.
\end{align*}
where $k, k_1, k_2\in C^{\infty}_c (G)$.

There are natural left and right actions of the commutative algebra
$C^\infty(M)$ on $C^{\infty}_c(G)$ given by the formulas
\[
a\cdot \sigma (\gamma)= a(r(\gamma))\sigma(\gamma), \quad \sigma
\cdot a(\gamma)= a(s(\gamma))\sigma(\gamma), \quad \gamma\in G,
\]
for any $a\in C^\infty(M)$ and $\sigma\in C^{\infty}_c(G)$.

We enlarge the algebra $C^{\infty}_c(G)$, introducing the unital
algebra
\[
\hat{C}^{\infty}(G)=C^{\infty}_c(G)+C^\infty(M)
\]
with the multiplication given by
\[
(k_1+a_1)(k_2+a_2)=k_1\ast k_2 + a_1\cdot k_2+ k_1\cdot a_2+a_1a_2,
\]
where $a_1a_2$ is the product of the functions $a_1$ and $a_2$.

We will also need the noncommutative analogue of the algebra of
differential forms on the leaf space of the foliation. Denote
$\Omega^j_\infty(G)=C^\infty_c(G, r^*\Lambda^jN^*\cF)$. There is a
product
\[\Omega^j_\infty(G) \times \Omega^k_\infty(G)
\stackrel{\wedge}{\to} \Omega^{j+k}_\infty(G)\] given, for any
$\omega\in \Omega^j_\infty(G)$ and $\omega_1\in \Omega^k_\infty(G)$,
by
\[
(\omega\wedge\omega_1)(\gamma)=\int_{G^y}\omega(\gamma_1)\wedge
H_{\gamma_1}[\omega_1(\gamma_1^{-1}\gamma)]d\nu^y(\gamma_1), \quad
\gamma\in G. \quad r(\gamma)=y.
\]
Here $H_{\gamma_1}: N^*_{s(\gamma_1)}\cF\to N^*_{s(\gamma_1)}\cF$ is
the linear holonomy map associated with $\gamma_1$.

One can also define natural left and right actions of the algebra
$\Omega^*_H(M)=C^\infty(M, \Lambda^*N^*\cF)$ of transverse
differential forms on $M$ on $\Omega^*_\infty(G)$ by the formulas
\[
a\wedge \omega (\gamma)= r^*a(\gamma)\wedge \omega(\gamma), \quad
\omega \wedge a(\gamma)= \omega(\gamma)\wedge
H_\gamma(s^*a(\gamma)), \quad \gamma\in G,
\]
for any $a\in \Omega^*_H(M)$ and $\omega\in \Omega^*_\infty(G)$.

We enlarge the algebra $\Omega^*_\infty(G)$, introducing the unital
algebra
\[
\hat{\Omega}^*_\infty(G)=\Omega^*_\infty(G)+\Omega^*_H(M),
\]
with the multiplication given by
\[
(\omega_1+a_1)\wedge (\omega_2+a_2)=\omega_1\wedge \omega_2 +
a_1\wedge \omega_2+ \omega_1 \wedge a_2 +a_1\wedge a_2,
\]
where $\omega_1\wedge \omega_2$ is the product of the forms
$\omega_1$ and $\omega_2$.

\subsection{Transverse differential}
Let $H\subset TM$ be a $q$-dimensional distribution such that
$TM=F\oplus H$. There is \cite{Co,Sau} the transverse
differentiation, which is a linear map
\[
D_H: \Omega_\infty^0(G)= C^\infty_c(G) \to
\Omega_\infty^1(G)=C^\infty_c(G, r^*N^*\cF),
\]
satisfying the condition
\[
D_H(k_1\ast k_2) = D_Hk_1\ast k_2+ k_1\ast D_Hk_2, \quad k_1, k_2
\in C^\infty_c(G).
\]
In this Section, we recall its definition.

The transverse distribution $H$ naturally defines a transverse
distribution $HG\cong r^*H$ on the foliated manifold $(G,\cG)$ and
the corresponding transversal de Rham differential $d_H:
C^\infty_c(G) \to C^\infty_c(G, r^*N^*\cF)$. For any $X\in H_y$,
there is a unique vector $\widehat{X}\in T_{\gamma}G$ such that $d
s(\widehat{X})= dh^{-1}_{\gamma}(X)$ and $d r(\widehat{X})=X$, where
$dh_{\gamma} : H_x\to H_y$ is the linear holonomy map associated
with $\gamma$. The space $H_\gamma G$ consists of all vectors of the
form $\widehat{X}\in T_{\gamma}G$ for different $X\in H_y$. In any
coordinate chart $W(\phi,\phi')$ on $G$, the distribution $H_\gamma
G$ consists of vectors $X\frac{\partial}{\partial
x}+X'\frac{\partial}{\partial x'}+Y\frac{\partial}{\partial y}$ such
that $X\frac{\partial}{\partial x}+Y\frac{\partial}{\partial y}\in
H_{(x,y)}$ and $X'\frac{\partial}{\partial
x'}+Y\frac{\partial}{\partial y}\in H_{(x',y)}$.

For any $f\in C^\infty_c(G)$, define $d_Hf\in C^\infty_c(G,
r^*N^*\cF)$ by
\[
d_Hf(X)=d f(\widehat{X}), \quad  X\in (r^*\tau)_\gamma \cong H_y,
\quad \gamma : x\to y,
\]
where $\widehat{X}\in H_\gamma G \subset T_{\gamma}G$ is a unique
vector such that $d s(\widehat{X})= dh^{-1}_{\gamma}(X)$ and $d
r(\widehat{X})=X$.

For the fixed smooth leafwise density $\alpha \in
C^\infty_c(M,|T\cF|)$, we define a transverse 1-form $k(\alpha)\in
C^\infty(M,H^*)\cong C^\infty(M, N^*\cF)$ as follows. Take an
arbitrary point $m \in M$ and $X\in H_{m}$. Let $\tilde{X}$ be an
arbitrary local projectable vector field, that coincides with $X$ at
$m$. In a foliated chart $\phi : U \to I^p\times I^q$ near $m$ such
that $\phi(m)=(x^0,y^0)$, one can write
\[
\tilde{X}(x,y)=\sum_{i=1}^pX^i(x,y)\frac{\partial}{\partial
x_i}+\sum_{j=1}^qY^j(y)\frac{\partial}{\partial y_j}.
\]
Then we put
\begin{multline*}
k(\alpha)(X)=\sum_{i=1}^pX^i(x^0,y^0)\frac{\partial f}{\partial
x_i}(x^0,y^0)+\sum_{j=1}^qY^j(y^0)\frac{\partial f}{\partial
x_j}(x^0,y^0)\\ +\sum_{i=1}^p \frac{\partial X^i}{\partial
x_i}(x^0,y^0) f(x^0,y^0).
\end{multline*}
It can be checked that this definition is independent of the choice
of a foliated chart $\phi$ and an extension $\tilde{X}$. If $M$ is
Riemannian, $\alpha $ is given by the induced leafwise Riemannian
volume form, and $H=F^\bot$, then $k(\alpha)$ coincides with the
mean curvature 1-form of $\cF$ (cf., for instance, \cite{Tondeur}).

For $f\in C^\infty_c(G)$, define $D_Hf\in C^\infty_c(G, r^*N^*\cF)$
as
\[
D_Hf(\gamma) = d_Hf(\gamma)  + \frac{1}{2}
(H_\gamma[s^*k(\alpha)(\gamma)] + r^*k(\alpha)(\gamma)) f(\gamma),
\quad \gamma\in G.
\]

Finally, note that the operator $D_H$ has a unique extension to a
differentiation of the differential graded algebra
$\Omega_\infty(G)$ (see \cite{Co,Sau}).

\section{Noncommutative Poisson geometry of foliations}
\subsection{Noncommutative Poisson structures}
Let $(M,\cF)$ be a transversely symplectic compact foliated
manifold, and $\omega$ the corresponding closed two-form of constant
rank on $M$. In this Section, we describe a Poisson structure on the
algebra $\hat{C}^{\infty}(G)$ associated to the foliation $\cF$
\cite{Block-Ge}.

First, we need some facts about connections on foliated manifolds.
Recall that there is a canonical flat connection
\[
\stackrel{\circ}{\nabla} : C^\infty(M,T\cF) \times
C^\infty(M,\tau)\to C^\infty(M,\tau)
\]
in the normal bundle $\tau$, defined along the leaves of $\cF$ (the
Bott connection). It is given by
\begin{equation}\label{e:Bott}
{\stackrel{\circ}\nabla}_X N=\theta(X)N=
P_\tau[X,\widetilde{N}],\quad X\in C^\infty(M,T\cF),\quad N\in
C^\infty(M,\tau),
\end{equation}
where $P_\tau : TM\to \tau$ is the natural projection and
$\widetilde{N}\in C^\infty(M, TM)$ is any vector field on $M$ such
that $P_\tau(\widetilde{N})= N$. Thus, the restriction of $\tau$ to
any leaf of $\cF$ is a flat vector bundle. The parallel transport in
$\tau$ along any leafwise path $\gamma:x\to y$ defined by
${\stackrel{\circ}\nabla}$ coincides with the linear holonomy map
$dh_{\gamma}:\tau_x\to\tau_y$.

A connection $\nabla : C^\infty(M,TM)\times C^\infty(M,\tau)\to
C^\infty(M,\tau)$ in the normal bundle $\tau$ is called adapted, if
its restriction to $C^\infty(M,T\cF)$ coincides with the Bott
connection ${\stackrel{\circ}\nabla}$.

One can construct an adapted connection, starting with an arbitrary
Riemannian metric $g_M$ on $M$. Denote by $\nabla^g$ the Levi-Civita
connection, defined by $g_M$. An adapted connection $\nabla$ is
given by
\begin{equation}\label{e:adapt}
\begin{aligned}
\nabla_XN &=P_\tau[X,\widetilde{N}],\quad X\in
C^\infty(M,T\cF),\quad N\in C^\infty(M,\tau)\\
\nabla_XN&=P_\tau\nabla^g_X\widetilde{N},\quad X\in
C^\infty(M,T\cF^{\bot}),\quad N\in C^\infty(M,\tau),
\end{aligned}
\end{equation}
where $\widetilde{N}\in C^\infty(M,TM)$ is any vector field such
that $P_\tau(\widetilde{N})=N$. One can show that the adapted
connection $\nabla$ described above is torsion-free.

An adapted connection $\nabla$ in the normal bundle $\tau$ is called
holonomy invariant, if, for any $X\in C^\infty(M,T\cF)$, $Y\in
C^\infty(M,TM)$, and $N\in C^\infty(M,\tau)$, we have
\[
(\theta(X)\nabla)_YN:=\theta(X)[\nabla_YN]-\nabla_{\theta(X)Y}N-
\nabla_Y[\theta_YN] =0.
\]
A holonomy invariant adapted connection in $\tau$ is called a basic
(or projectable) connection.

If the foliation $\cF$ is Riemannian and $g_M$ is a bundle-like
metric, then the connection $\nabla$ defined by \eqref{e:adapt} is a
basic connection. There are topological obstructions for the
existence of basic connections for an arbitrary foliations found
independently by Kamber-Tondeur and Molino.

We will assume the following:

\begin{hyp}\label{hyp}
There exists a basic connection on the normal bundle $\tau$ of
$\cF$.
\end{hyp}

The two-form $\omega$ induces an isomorphism $\#_\omega$ between the
bundle $\tau$ and $\tau^*$:
\[
\langle \#_\omega X, Y\rangle=\omega(\widetilde{X},\widetilde{Y}),
\quad X, Y\in \tau,
\]
where $\widetilde{X}\in TM$ and $\widetilde{Y}\in TM$ are such that
$P_\tau(\widetilde{X})=X$. Thus, we have a skew-symmetric form on
$\tau^*$, which we denote by $\Lambda$:
\[
\Lambda(\#_\omega X,\#_\omega Y)=\omega(X,Y), \quad X,Y\in \tau.
\]
It is shown in \cite{Block-Ge} that, under Hypotheses~\ref{hyp}, $M$
has a presymplectic connection, that is, a basic connection $\nabla$
on $\tau$ such that $\nabla\Lambda=0$. From now on, we will assume
that $\nabla$ is a presymplectic connection.

The definition of a noncommutative Poisson structure depends on a
choice  of a $q$-dimensional distribution $H\subset TM$ such that
$TM=F\oplus H$. The Poisson bracket of two functions $k_1, k_2\in
C^\infty_c(G)$ is defined by the formula
\[
\Pi_H(k_1,k_2)=\Lambda(D_Hk_1, D_Hk_2)
\]
or
\[
\Pi_H(k_1,k_2)(\gamma)=\int_{G^y}\langle \Lambda_{y},
D_Hk_1(\gamma_1)\wedge H_{\gamma_1}
[D_Hk_2(\gamma_1^{-1}\gamma)]\rangle d\nu^y(\gamma_1),\quad
\gamma\in G^y.
\]

By \cite{Block-Ge}, this Poisson bracket $\Pi_H$ satisfies
(\ref{P1}) and (\ref{P2}) with $\Pi_1$ defined as follows. Let $
\nabla : C^\infty(M,N^*{\cF}) \rightarrow C^\infty_c(M,N^*\cF\otimes
N^*\cF) $ be the covariant derivative on the bundle $N^*\cF$ defined
by the connection $\nabla$. It gives rise to an operator
\[
\nabla:C^\infty_c(G,r^*N^*{\cF}) \rightarrow
C^\infty_c(G,r^*N^*\cF\otimes r^*N^*\cF).
\]
Denote by
\[
D^2=\nabla \circ D_H :C^\infty_c(G)\rightarrow
C^\infty_c(G,r^*N^*\cF\otimes r^*N^*\cF)
\]
the composition of $D_H:C^\infty_c(G)\rightarrow
C^\infty_c(G,r^*N^*\cF)$ with $\nabla$; $D^2$ takes values in
$C^\infty_c(G,S^2r^*N^*\cF)$ since $\nabla$ is torsion-free. Then
$\Pi_1$ is a two-chain on $C^\infty_c(G)$ defined by the formula
\[
\Pi_1(k_1,k_2)=\Lambda\otimes\Lambda(D^2k_1\ast D^2k_2),\quad k_1,
k_2\in C^\infty_c(G).
\]

We extend $\Pi_H$ to the algebra $\hat{C}^{\infty}(G)$ by the
formula
\begin{multline}\label{e:Pi}
\Pi_H(k_1+a_1, k_2+a_2)\\ =\Pi_H(k_1, k_2) + \Lambda(d_Ha_1, D_H
k_2)+ \Lambda(D_Hk_1, d_Ha_2)+\Lambda(d_Ha_1, d_Ha_2),
\end{multline}
where, for any $\gamma\in G$,
\begin{gather*}
\Lambda(d_Ha_1, D_H k_2)(\gamma)=\Lambda_y(d_Ha_1(r(\gamma)), D_H
k_2(\gamma)),\\
\Lambda(D_Hk_1, d_Ha_2)(\gamma)=\Lambda_y(D_Hk_1(\gamma),
H_\gamma[d_Ha_2(s(\gamma))]).
\end{gather*}

 It is easy to see that $\Pi_H$ is a noncommutative Poisson
structure on $\hat{C}^{\infty}(G)$ in the sense of
Definition~\ref{d:Poisson}.

\subsection{Transverse Hamiltonian flows}
As above, we suppose that $(M,\cF)$ is a transversely symplectic
compact foliated manifold and $\omega$ is the corresponding closed
two-form of constant rank on $M$. A Hamiltonian on the singular
symplectic manifold $M/\cF$ is given by a $C^\infty$ function $h$ on
$M$, which is constant on each leaf of the foliation $\cF$. It is
easy to see that $h$ belongs to the center of the algebra
$\hat{C}^{\infty}(G)$. The purpose of this section is to give an
explicit geometric description of the associated Hamiltonian
derivation $X_h$.

First, we recall (see \cite{Gotay}) that for any presymplectic
manifold $(M,\omega)$ there is a symplectic manifold $(\Phi,\eta)$
and an embedding $i:M\to\Phi$ such that $\omega=i^*\eta$ and $M$ is
a coisotropic submanifold of $\Phi$. Moreover, such a coisotropic
embedding is unique up to local symplectomorphism about $M$. Its
construction makes use of an auxiliary choice of a distribution
$H\subset TM$ such that $TM=H \oplus T\cF$. Such a distribution
yields an embedding $j$ of $T^*\cF$ in $T^*M$. Let $\pi : T^*\cF\to
M$ be the natural projection. Let $j^*\omega_{T^*M}$ be the
pull-back of the canonical symplectic form $\omega_{T^*M}$ on $T^*M$
to $T^*\cF$. Then one can take
\[
\eta=\pi^*\omega+j^*\omega_{T^*M}.
\]
It is easy to see that the restriction of $\eta$ to $M$ equals
$\omega$. The manifold $\Phi$ is defined to be a tubular
neighborhood of the zero section $M$ in $T^*\cF$ so that $\eta$
restricted to $\Phi$ is non-degenerate.

Remark that the restricted tangent bundle $T_M\Phi$ has the
canonical decomposition
\[
T_M\Phi\cong TM \oplus T^*\cF.
\]
Thus, we have
\[
T_M\Phi\cong H \oplus T\cF\oplus T^*\cF.
\]
Denote by $p$ the induced projection $T_M\Phi\to  T\cF\oplus
T^*\cF$. For $m\in M$, let $\omega_F$ denote the canonical
symplectic structure on $T\cF\oplus T^*\cF$:
\[
\omega_F(f_1\oplus f^*_1, f_2\oplus f^*_2)=\langle f^*_2, f_1\rangle
-\langle f^*_1, f_2\rangle.
\]
Then the restriction of $\eta$ to $T_M\Phi$ is described as
\[
\eta=\pi^*\omega+\omega_{F}\circ (p\times p).
\]
Thus, for any $X=\pi_*(X)+f^*_X\in T_M\Phi$ and $Y=\pi_*(Y)+f^*_Y\in
T_M\Phi$, we have
\[
\eta(X,Y)=\omega (\pi_*(X), \pi_*(Y))+\langle f^*_Y,
p_F(\pi_*(X))\rangle -\langle f^*_X, p_F(\pi_*(Y))\rangle.
\]

Given a $C^\infty$ function $h$ on $M$, which is constant on each
leaf of $\cF$, we extend it to a smooth function $\tilde{h}$ on
$\Phi$. Let $v_{\tilde{h}}$ be the Hamiltonian vector field of the
function $\tilde{h}$ on $\Phi$. Recall that by definition we have
\[
i_{v_{\tilde{h}}}\eta=d\tilde{h}.
\]
Then (see, for instance, \cite{LM87}) the submanifold $M$ of $\Phi$
is invariant under the flow of $v_{\tilde{h}}$. Indeed, for
$Y=\pi_*(Y)\in TM\subset  T_M\Phi, Y=p_F(Y)$, we have
\[
0=d\tilde{h}(Y)=\eta(v_{\tilde{h}},Y)=-\langle
v_{\tilde{h}}-\pi_*(v_{\tilde{h}}), p_F(\pi_*(Y))\rangle.
\]
Therefore, we conclude that $v_{\tilde{h}}=\pi_*(v_{\tilde{h}})$.

It is easy to see that $v_h$ depends only on $h$ and $d\tilde{h}$
restricted to $T^*\cF\subset T_M\Phi $. For any $Y=f^*_Y\in
T^*\cF\subset T_M\Phi$, we have
\[
d\tilde{h}(Y)=\eta(v_{\tilde{h}},Y)=\langle f^*_Y,
p_F(v_{\tilde{h}})\rangle.
\]
Thus, we see that
\[
d\tilde{h}\left|_{T^*\cF}\right.=p_F(v_{\tilde{h}})\in
C^\infty(M,T\cF).
\]

Finally, if we denote by $v_h$ the restriction of $v_{\tilde{h}}$ to
$M$, then one can show that the flow of $v_h$ on $M$ preserves the
foliation $\cF$, that is, it takes a leaf of $\cF$ into a leaf.
Therefore, there is a natural lift of $v_h$ to a vector field
$\hat{v}_h$ on $G$ such that for any $\gamma\in G$,
$s_{\ast}(\hat{v}_h(\gamma))=v_h(s(\gamma))$ and
$r_{*}(\hat{v}_h(\gamma))=dh_{\gamma}[v_h(s(\gamma))]=v_h(r(\gamma))$,
where $dh_{\gamma}$ is the differential of the holonomy map along
$\gamma$. In local foliated coordinates, $v_h$ has a form
\[
v_h(x,y) = \sum_{j=1}^p X^j(x,y) \frac{\partial }{\partial x_j} +
\sum_{k=1}^q Y^k(y)\frac{\partial }{\partial y_k}, \quad (x,y) \in
I^p\times I^q.
\]
and $\hat{v}_h$ is given by
\begin{multline*}
\hat{v}_h(x,x',y) = \sum_{j=1}^p X^j(x,y) \frac{\partial }{\partial
x_j} +\sum_{j=1}^p X^j(x^\prime,y) \frac{\partial }{\partial
x^\prime_j}+ \sum_{k=1}^q Y^k(y)\frac{\partial }{\partial y_k},\\
\quad (x,x',y) \in I^p\times I^p\times I^q.
\end{multline*}

Define an operator ${\mathcal L}_{\hat{v}_h}$ on the space
$C^\infty_c(G)$ by the formula
\[
{\mathcal L}_{\hat{v}_h}f = \hat{v}_h f  + \frac{1}{2}
\left(s^*[k(\alpha)(\hat{v}_h)] + r^*[k(\alpha)(\hat{v}_h)]\right)
f, \quad f\in C^\infty_c(G).
\]
It coincides with the Lie derivative by $\hat{v}_h$ on the space
$C^{\infty}_{c}(G, |T{\mathcal G}|^{1/2})$ of leafwise
half-densities on the holonomy groupoid $G$. In a foliated chart,
for any $k\in C^{\infty}_{c}(G)$, we have
\begin{multline*}
\cL_{\hat{v}_h}(k)\\ = \Big(\hat{v}_h
k(x,x',y)+\frac{1}{2}\sum_{j=1}^p \frac{\partial X^j}{\partial
x_j}(x,y) k(x,x^\prime,y) +\frac{1}{2}\sum_{j=1}^p \frac{\partial
X^j}{\partial x^\prime_j} (x^\prime,y)
k(x,x^\prime,y)\Big),\\
\quad (x,x',y) \in I^p\times I^p\times I^q.
\end{multline*}

We arrive at the main result of the paper.

\begin{thm}
Suppose that $(M,\cF)$ is a transversely symplectic compact foliated
manifold such that there exists a basic connection on the normal
bundle $\tau$ of $\cF$. Let $H$ be a $q$-dimensional distribution on
$M$ such that $TM=F\oplus H$, $\Pi_H$ the noncommutative Poisson
structure on the algebra $\hat{C}^{\infty}(G)$ defined
by~\eqref{e:Pi}, and $i_H : M\to\Phi_H$ the corresponding
coisotropic embedding into a symplectic manifold $(\Phi_H,\eta_H)$.

Let $h$ be a $C^\infty$ function on $M$, which is constant on each
leaf of $\cF$, and $\tilde{h}$ its extension to a smooth function on
$\Phi_H$ such that $d\tilde{h}\left|_{T^*\cF}\right.=0$.

The Hamiltonian derivation $X_h$ of the algebra $\hat{C}^\infty(G)$
associated to the Hamiltonian $h$ and the noncommutative Poisson
structure $\Pi_H$ coincides with ${\mathcal L}_{\hat{v}_h}+v_h$:
\[
X_h(k+a)={\mathcal L}_{\hat{v}_h}(k)+v_h(a), \quad k\in
C^\infty_c(G),\quad a\in C^\infty(M).
\]
\end{thm}

\begin{proof}
The key step in the proof is the following lemma.

\begin{lem}
For any $a\in C^\infty(M)$, we have
\[
\Lambda(d_Hh,d_Ha)=v_h(a).
\]
\end{lem}

\begin{proof}
Denote by $\#_\eta : T\Phi\to T^*\Phi$ the isomorphism induced by
the two-form $\eta$:
\[
\langle \#_\eta X, Y\rangle=\eta(X,Y), \quad X, Y\in T\Phi.
\]
and by $\Lambda_\Phi$ the induced two-form on $T^*\Phi$
\[
\Lambda_\Phi(\#_\eta X,\#_\eta Y)=\eta(X,Y), \quad X,Y\in T\Phi.
\]
It is easy to see that $\#_\eta$ maps $TM$ to $N^*\cF$, and the
kernel of the map $\#_\eta : TM \to N^*\cF$ coincides with $T\cF$.
Thus, we have the induced map $\bar{\#}_\eta : TM/TF \to N^*\cF$,
which is equal to $\#_\omega$.

Using these facts, we easily derive that
\[
\Lambda_\Phi(\nu_1,\nu_2)=\Lambda(\nu_1,\nu_2), \quad \nu_1,\nu_2\in
N^*\cF.
\]

On the other hand, for a given function $a\in C^\infty(M)$ take its
extension $\tilde{a}\in C^\infty(\Phi)$ to $\Phi$ such that $
d\tilde{a}\left|_{T^*\cF}\right.=0$. Then by definition we have
\[
\Lambda_\Phi(d\tilde{h},d\tilde{a})=d\tilde{a}(v_{\tilde{h}}).
\]
Let us restrict both sides of this identity to $M$. Then by
assumption the restriction of $d\tilde{h}\in C^\infty(\Phi,
T^*\Phi)$ to $M$ coincides with $dh\in C^\infty(M, T^*M)\subset
C^\infty(M, T^*_M\Phi)$. Moreover, we have $dh=d_Hh$. Similarly, the
restriction of $d\tilde{a}\in C^\infty(\Phi, T^*\Phi)$ to $M$
coincides with $da\in C^\infty(M, T^*M)\subset C^\infty(M,
T^*_M\Phi)$. Since $v_{\tilde{h}}\left|_{M}\right.=v_h\in
C^\infty(M, TM)\subset C^\infty(M, T_M\Phi)$, in particular, this
implies that $d\tilde{a}(v_{\tilde{h}})\left|_{M}\right.=da(v_h)$.
We arrive at the identity
\[
\Lambda_\Phi(d_Hh,da)=da(v_h).
\]
It remains to show that
\[
\Lambda_\Phi(d_Hh,d_Fa)=0.
\]
Given $X\in T_M\cF$, $Y=\pi_*(Y)+f^*_Y\in T_M\Phi$, we have
\[
\langle \#_\eta X, Y\rangle=\eta(X,Y)=\langle f^*_Y, X\rangle..
\]
Therefore, $\#_\eta X\in T^*\cF$, and $\#_\eta^{-1} : T^*\cF\to
T\cF$. So $\#_\eta^{-1}d_Hh\in TM$, $\#_\eta^{-1}d_Fa\in T\cF$, and
we obtain
\[
\Lambda_\Phi(d_Hh,d_Fa)=\eta(\#_\eta^{-1}d_Hh,\#_\eta^{-1}d_Fa)
=\omega(\#_\eta^{-1}d_Hh,\#_\eta^{-1}d_Fa)=0,
\]
that completes the proof of the lemma.
\end{proof}

By this lemma, it follows easily that, for any  $a\in C^\infty(M)$,
\[
\frac12(\Pi(h,a)-\Pi(a,h))=v_h(a)
\]
and, for any $k\in C^\infty_c(G)$,
\[
\frac12(\Pi(h,k)-\Pi(k,h))={\mathcal L}_{\hat{v}_h}(k),
\]
that completes the proof.
\end{proof}


\begin{thebibliography}{10}

\bibitem{Block-Ge}
J. Block and E. Getzler. {\em Quantization of foliations}. In
Proceedings of the {XX}th {I}nternational {C}onference on
  {D}ifferential {G}eometric {M}ethods in {T}heoretical {P}hysics, {V}ol.\ 1, 2
  ({N}ew {Y}ork, 1991), 471--487. World Sci. Publ., River Edge, NJ,
  1992.

\bibitem{Co79}
A. Connes. {\em Sur la th\'eorie non commutative de
l'int\'egration.} In {Alg\`ebres d'op\'erateurs ({S}\'em., {L}es
{P}lans-sur-{B}ex, 1978)}, 19--143. {Lecture Notes in Math.}, Vol.
725. Springer, Berlin, 1979.


\bibitem{Co}
A. Connes. {\em Noncommutative geometry}. Academic Press Inc., San
Diego, CA, 1994.

\bibitem{Gersten63}
M.~Gerstenhaber. {\em The cohomology structure of an associative
ring.} {Ann. of Math.} \textbf{78} (1963), 267--288.

\bibitem{Ginzburg}
V. Ginzburg. {\em Non-commutative symplectic geometry, quiver
varieties, and operads.} {Math. Res. Lett.} \textbf{8} (2001),
377--400.

\bibitem{Ginzburg:lectures}
V. Ginzburg. {\em Lectures on noncommutative geometry.} Preprint
arXiv:math.AG/0506603, 2005.

\bibitem{Gotay}
M.~J. Gotay. {\em On coisotropic imbeddings of presymplectic
manifolds.} {Proc. Amer. Math. Soc.} \textbf{84} (1982), 111--114.

\bibitem{Halbout-Tang}
G. Halbout and X. Tang. {\em Noncommutative {Poisson} structures on
orbifolds.} Preprint arXiv:math/0606436, 2006.

\bibitem{kontsevich93}
M. Kontsevich. {\em Formal (non)commutative symplectic geometry.} In
{The {G}el'fand {M}athematical {S}eminars, 1990--1992}, 173--187.
Birkh\"auser Boston, Boston, MA, 1993.

\bibitem{dgtrace-eng}
Yu.~A. Kordyukov. {\em The trace formula for transversally elliptic
operators on {R}iemannian foliations.} {Algebra i Analiz}
\textbf{12} (2000), No. 3, 81--105.

\bibitem{egorgeo}
Yu.~A. Kordyukov. {\em Egorov's theorem for transversally elliptic
operators on foliated manifolds and noncommutative geodesic flow.}
{Math. Phys. Anal. Geom.} \textbf{8} (2005), No. 2, 97--119.

\bibitem{matrix-egorov}
Yu~A. Kordyukov. {\em The {E}gorov theorem for transverse
{D}irac-type operators on foliated manifolds.} {J. Geom. Phys.},
\textbf{57} (2007), 2345--2364.

\bibitem{survey}
Yu~A. Kordyukov. {\em Noncommutative geometry of foliations.} J.
K-Theory \textbf{2} (2008), 219--327.

\bibitem{LM87}
P. Libermann and Ch.-M. Marle.
\newblock {\em Symplectic geometry and analytical mechanics.}
D. Reidel Publishing Co., Dordrecht, 1987.

\bibitem{Sau}
J.-L. Sauvageot. {\em Semi-groupe de la chaleur transverse sur la
${C}^*$-algebre d'un feuilletage riemannien.} {J. Funct. Anal.}
\textbf{142} (1996), 511--538.

\bibitem{tang:gafa}
X. Tang. {\em Deformation quantization of pseudo-symplectic
({P}oisson) groupoids.} {Geom. Funct. Anal.} \textbf{16} (2006),
731--766.

\bibitem{Tondeur}
Ph. Tondeur. {\em Geometry of foliations.} Birkh\"auser, Basel,
1997.

\bibitem{Vaisman}
I. Vaisman. {\em Lectures on the geometry of Poisson manifolds.}
Birkh\"auser, Basel, 1994.

\bibitem{Vey75}
J. Vey. {\em D\'eformation du crochet de {P}oisson sur une
vari\'et\'e symplectique.} {Comment. Math. Helv.} \textbf{50}
(1975), 421--454.

\bibitem{Xu}
P. Xu. {\em Noncommutative {P}oisson algebras.} {Amer. J. Math.}
\textbf{116} (1994), 101--125.

\end{thebibliography}

\end{document}